# Some elementary observations on Narayana polynomials and related topics II: q-Narayana polynomials

Johann Cigler

Fakultät für Mathematik

Universität Wien

johann.cigler@univie.ac.at

**Abstract**

We show that $q$-Catalan numbers, $q$-central binomial coefficients and $q$-Narayana polynomials are moments of $q$-analogues of Fibonacci and Lucas polynomials and related polynomials.

## 1. Introduction

This note is a supplement to part I ([4]). Let

$$F_n(x) = \sum_{k=0}^{\left\lfloor \frac{n}{2} \right\rfloor} (-1)^k \binom{n-k}{k} x^{n-2k}, \qquad (1.1)$$

$n \in \mathbb{N}$, be Fibonacci polynomials and define a linear functional $L$ by

$$L(F_n(x)) = [n = 0]. \qquad (1.2)$$

Then the moments $L(x^{2n})$ are the Catalan numbers

$$L(x^{2n}) = C_n = \frac{1}{n+1}\binom{2n}{n}. \qquad (1.3)$$

This well-known fact is the special case $t = 1$ of the following result (cf. [4] and the literature cited there):

The Narayana polynomials

$$C_n(t) = \sum_{k \geq 0} \binom{n}{k}\binom{n}{k+1}\frac{1}{n}t^k \qquad (1.4)$$

can be represented as moments $L(x^{2n})$ of the linear functional $L$ defined by $L(F_n(x,t)) = [n = 0],$ where

$$F_n(x,t) = \sum_{k=0}^{\left\lfloor \frac{n}{2} \right\rfloor} (-1)^k \sum_{j=0}^{k} \binom{\left\lfloor \frac{n}{2} \right\rfloor - j}{k-j}\binom{\left\lfloor \frac{n-1}{2} \right\rfloor - k + j}{j} t^j x^{n-2k} \qquad (1.5)$$

are generalized Fibonacci polynomials which satisfy the recurrence

$$F_{2n}(x,t) = xF_{2n-1}(x,t) - F_{2n-2}(x,t),$$
$$F_{2n+1}(x,t) = xF_{2n}(x,t) - tF_{2n-1}(x,t).$$
(1.6)

One of our purposes is to give a nice $q-$analogue of this result. More precisely we define nice polynomials $F_n(x,t,q)$ such that the linear functional defined by $L(F_n(x,t,q)) = [n=0]$ has as moments the $q-$Narayana polynomials

$$L(x^{2n}) = C_n(t,q) = \sum_{k=0}^{n} \begin{bmatrix} n \\ j \end{bmatrix} \begin{bmatrix} n \\ j+1 \end{bmatrix} \frac{1}{[n]} q^{j^2+j} t^j,$$
(1.7)

which for $t=1$ reduce to the $q-$Catalan numbers

$$C_n(q) = C_n(1,q) = \frac{1}{[n+1]} \begin{bmatrix} 2n \\ n \end{bmatrix}.$$
(1.8)

We will always suppose that $0 < q < 1$ and use the notations $[n] = [n]_q = \frac{1-q^n}{1-q}$ and

$$\begin{bmatrix} n \\ k \end{bmatrix} = \begin{bmatrix} n \\ k \end{bmatrix}_q = \frac{(q;q)_n}{(q;q)_k (q;q)_{n-k}} \text{ for } 0 \le k \le n.$$

It is well known that for each sequence $(a_n)_{n \ge 0}$ with $a_0 = 1$ such that all Hankel determinants $\det(a_{i+j})_{i,j=0}^n \ne 0$ there are uniquely defined monic polynomials $p_n(x)$ of degree $n$, $n \in \mathbb{N}$, which are orthogonal with respect to the linear functional $L$ defined by $L(x^n) = a_n$. Orthogonality means that $L(p_n p_m) = 0$ for $n \ne m$ and $L(p_n^2) \ne 0$. Since $p_0(x) = 1$ this implies that $L(p_n) = [n=0]$.

Moreover there are uniquely determined numbers $a(n,k)$ such that

$$x^n = \sum_{k=0}^{n} a(n,k) p_k(x).$$
(1.9)

By applying $L$ we see that $a(n,0) = L(x^n) = a_n$.

For $q=1$ the polynomials $F_n(x,t)$ are the orthogonal polynomials with moments $L(x^{2n}) = C_n(x,t)$ and $L(x^{2n+1}) = 0$.

Let now $L$ be the linear functional defined by $L(x^{2n}) = C_n(t,q)$ and $L(x^{2n+1}) = 0$.

It would be tempting to consider also in this case the corresponding orthogonal polynomials, but it seems that there is no simple explicit formula or recurrence for them, not even for $t=1$.

In that case the first orthogonal polynomials are

$$1, \quad x, \quad x^2-1, \quad x^3-(1+q^2)x, \quad x^4-(q+q^2+q^4)x^2-1+q+q^4,$$

$$x^5 - \frac{(1+q)(1-q+q^2)(-1+q^3+q^5+q^7)}{-1+q+q^4}x^3 + \frac{q(1-q+q^2)(-1-q+q^5+q^6+q^7+q^8+q^9)}{-1+q+q^4}x.$$

The first numbers $a(n,k)$ are

```
1
0           1
1           0                            1
0           1 + q²                       0                   1
1 + q²      0                            q (1 + q + q³)      0                                          1
0           (1 - q + q²) (1 + q + q² + q³ + q⁴)    0         (1+q) (1-q+q²) (-1+q³+q⁵+q⁷)               0    1
                                                             ─────────────────────────────
                                                                     -1+q+q⁴
```

Fortunately there also exist "nice" polynomials with the same moments. Let us consider first the case $q=1$. The orthogonal polynomials for the linear functional $L$ defined by

$L(x^{2n}) = C_n$ and $L(x^{2n+1}) = 0$ are the Fibonacci polynomials $F_n(x) = \sum_{k=0}^{n}(-1)^k \binom{n-k}{k} x^{n-2k}$.

They satisfy $x^n = \sum_{k=0}^{n} a(n,k) F_k(x)$ with $a(2n+k,k) = \binom{2n+k}{n} - \binom{2n+k}{n-1} = \frac{k+1}{n+k+1}\binom{2n+k}{n}$

and $a(n,k) = 0$ else.

Let me sketch how to find a nice $q$-analogue of this situation. It is easier to begin with $a(n,k)$. A natural $q$-analogue is

$$a(2n+k,k,q) = \frac{1}{q^n}\left(\begin{bmatrix}2n+k\\n\end{bmatrix} - \begin{bmatrix}2n+k\\n-1\end{bmatrix}\right) = \frac{[k+1]}{[n+k+1]}\begin{bmatrix}2n+k\\n\end{bmatrix} \text{ and } a(n,k,q) = 0 \text{ else. The}$$

first terms are

```
1
0           1
1           0                            1
0           1 + q                        0                   1
1 + q²      0                            1 + q + q²          0                   1
0           1 + q + q² + q³ + q⁴         0                   (1 + q) (1 + q²)    0    1
```

Now we are looking for the polynomials $F_n(x,q)$ such that

$$x^n = \sum_{k=0}^{n} a(n,k,q) F_k(x,q). \tag{1.10}$$

Their coefficients are given by the inverse matrix $(a(i,j,q))^{-1}$. Fortunately this also turns out to be nice. The first terms are

$$\begin{pmatrix} 1 & 0 & 0 & 0 & 0 & 0 \\ 0 & 1 & 0 & 0 & 0 & 0 \\ -1 & 0 & 1 & 0 & 0 & 0 \\ 0 & -1-q & 0 & 1 & 0 & 0 \\ q & 0 & -1-q-q^2 & 0 & 1 & 0 \\ 0 & q(1+q+q^2) & 0 & -(1+q)(1+q^2) & 0 & 1 \end{pmatrix}$$

Thus the sequence $(F_n(x,q))_{n \geq 0}$ begins with

$1,\ x,\ x^2-1,\ x^3-(1+q)x,\ x^4-(1+q+q^2)x^2+q,$
$x^5-(1+q)(1+q^2)x^3+q(1+q+q^2)x,\cdots.$

$F_n(x,q)$ can of course also be computed inductively by $F_n(x,q) = x^n - \sum_{k=1}^{n} a(n,k,q)F_k(x,q)$.

It is now easy to guess that in general

$$F_n(x,q) = \sum_{j=0}^{\left\lfloor \frac{n}{2} \right\rfloor} (-1)^j q^{\binom{j}{2}} \begin{bmatrix} n-j \\ j \end{bmatrix} x^{n-2j}. \tag{1.11}$$

**Remark**

Note that these $q-$Fibonacci polynomials, which have been considered in [2] and [3], are not orthogonal. There are also nice orthogonal $q-$analogues of the Fibonacci polynomials, i.e. the Carlitz $q-$Fibonacci polynomials $\sum_{j=0}^{\left\lfloor \frac{n}{2} \right\rfloor} (-1)^j q^{j^2} \begin{bmatrix} n-j \\ j \end{bmatrix} x^{n-2j}$. Their moments are the Carlitz $q-$Catalan numbers, but unfortunately these have no closed formula.

In this note we first recall some results about the above mentioned class of non-orthogonal $q-$Fibonacci and $q-$Lucas polynomials whose moments are $q-$Catalan numbers and $q-$central binomial coefficients and then propose "nice" $q-$analogues of the generalized Fibonacci and Lucas polynomials of part I such that the corresponding moments are $q-$Narayana polynomials $C_n(t,q) = \sum_{k=0}^{n} \begin{bmatrix} n \\ j \end{bmatrix} \begin{bmatrix} n \\ j+1 \end{bmatrix} \frac{1}{[n]} q^{j^2+j} t^j$ and $M_n(t,q) = \sum_{j=0}^{n} q^{j^2} \begin{bmatrix} n \\ j \end{bmatrix}^2 t^j$.

## 2. Some background material

Let us first state some known results (cf. [1],[2],[3]). As already mentioned the $q$ – Fibonacci polynomials

$$F_n(x,q) = \sum_{j=0}^{\lfloor \frac{n}{2} \rfloor} (-1)^j q^{\binom{j}{2}} \begin{bmatrix} n-j \\ j \end{bmatrix} x^{n-2j} \tag{2.1}$$

satisfy

$$x^n = \sum_{k=0}^{\lfloor \frac{n}{2} \rfloor} \frac{1}{q^k} \left( \begin{bmatrix} n \\ k \end{bmatrix} - \begin{bmatrix} n \\ k-1 \end{bmatrix} \right) F_{n-2k}(x,q). \tag{2.2}$$

If we define a linear functional $L$ by $L(F_n(x,q)) = [n=0]$ then we get

$$L(x^{2n}) = C_n(q) = \begin{bmatrix} 2n \\ n \end{bmatrix} \frac{1}{[n+1]}, \tag{2.3}$$

where $C_n(q) = \frac{1}{q^n}\left(\begin{bmatrix} 2n \\ n \end{bmatrix} - \begin{bmatrix} 2n \\ n-1 \end{bmatrix}\right) = \begin{bmatrix} 2n \\ n \end{bmatrix}\frac{1}{[n+1]}$ is a $q$ – analogue of the Catalan numbers $C_n = \binom{2n}{n}\frac{1}{n+1}$.

The $q$ – Lucas polynomials

$$L_n(x,q) = \sum_{k=0}^{\lfloor \frac{n}{2} \rfloor} (-1)^k q^{\binom{k}{2}} \frac{[n]}{[n-k]} \begin{bmatrix} n-k \\ k \end{bmatrix} x^{n-2k} \tag{2.4}$$

for $n > 0$ and $L_0(x,q) = 1$ satisfy

$$x^n = \sum_{k=0}^{\lfloor \frac{n}{2} \rfloor} \begin{bmatrix} n \\ k \end{bmatrix} L_{n-2k}(x,q). \tag{2.5}$$

If we define the linear functional $M$ by $M(L_n(x,q)) = [n=0]$ then we get

$$M(x^{2n}) = M_n(q) = \begin{bmatrix} 2n \\ n \end{bmatrix} \tag{2.6}$$

is a central $q$ – binomial coefficient.

It is easy to verify that

$$L_n(x,q) = F_n(x,q) - q^{n-1} F_{n-2}(x,q) \tag{2.7}$$

for $n > 1$. Moreover $L_0(x,q) = F_0(x,q) = 1$ and $L_1(x,q) = F_1(x,q) = x$.

These results can be proved with an inversion formula by L. Carlitz [1]. Another proof is in [2] and [3]. Carlitz uses the fact that

$$c(n,k) := \sum_{j=0}^{\min(k,n-k)} (-1)^j q^{\binom{j+1}{2}} \begin{bmatrix} k \\ j \end{bmatrix} \begin{bmatrix} n-j \\ k \end{bmatrix} = 1 \tag{2.8}$$

for $0 \le k \le n$.

To prove this let $U$ be the linear operator on $\mathbb{C}[q^x]$ defined by $U\begin{bmatrix} x \\ k \end{bmatrix} = \begin{bmatrix} x-1 \\ k \end{bmatrix}$ for integers $k > 0$ and $U1 = 0$.

Then $c(n,k) := \sum_{j=0}^{k} (-1)^j q^{\binom{j+1}{2}} \begin{bmatrix} k \\ j \end{bmatrix} U^j \begin{bmatrix} n \\ k \end{bmatrix} = (1-qU)(1-q^2U)\cdots(1-q^kU) \begin{bmatrix} n \\ k \end{bmatrix} = \begin{bmatrix} n-k \\ 0 \end{bmatrix} = 1$

because $(1-q^kU)\begin{bmatrix} n \\ k \end{bmatrix} = \begin{bmatrix} n \\ k \end{bmatrix} - q^k \begin{bmatrix} n-1 \\ k \end{bmatrix} = \begin{bmatrix} n-1 \\ k-1 \end{bmatrix}$.

**Lemma (Carlitz [1], Theorem 7)**

If $u(n) = \sum_{k=0}^{\lfloor n/2 \rfloor} \frac{1}{q^k}\left(\begin{bmatrix} n \\ k \end{bmatrix} - \begin{bmatrix} n \\ k-1 \end{bmatrix}\right) v(n-2k)$

*then*

$v(n) = \sum_{j=0}^{\lfloor n/2 \rfloor} (-1)^j q^{\binom{j}{2}} \begin{bmatrix} n-j \\ j \end{bmatrix} u(n-2j).$

**Proof**

$\sum_{j=0}^{\lfloor n/2 \rfloor} (-1)^j q^{\binom{j}{2}} \begin{bmatrix} n-j \\ j \end{bmatrix} u(n-2j) = \sum_{j=0}^{\lfloor n/2 \rfloor} (-1)^j q^{\binom{j}{2}} \begin{bmatrix} n-j \\ j \end{bmatrix} \sum_{\ell=0}^{\lfloor (n-2j)/2 \rfloor} \frac{1}{q^\ell}\left(\begin{bmatrix} n-2j \\ \ell \end{bmatrix} - \begin{bmatrix} n-2j \\ \ell-1 \end{bmatrix}\right) v(n-2j-2\ell)$

$= \sum_{k=0}^{\lfloor n/2 \rfloor} v(n-2k) \sum_{j=0}^{k} (-1)^j q^{\binom{j}{2}} \begin{bmatrix} n-j \\ j \end{bmatrix} \frac{1}{q^{k-j}} \left(\begin{bmatrix} n-2j \\ k-j \end{bmatrix} - \begin{bmatrix} n-2j \\ k-j-1 \end{bmatrix}\right)$

$= v(n) + \sum_{k=1}^{\lfloor n/2 \rfloor} v(n-2k) \frac{1}{q^k} \left( \sum_{j=0}^{k} (-1)^j q^{\binom{j+1}{2}} \begin{bmatrix} k \\ j \end{bmatrix} \begin{bmatrix} n-j \\ k \end{bmatrix} - \sum_{j=0}^{k} (-1)^j q^{\binom{j+1}{2}} \begin{bmatrix} k-1 \\ j \end{bmatrix} \begin{bmatrix} n-j \\ k-1 \end{bmatrix} \right) = v(n)$

by (2.8).

If we choose $u(n) = x^n$ we get (2.1).

By (2.7) and

$$q^k \begin{bmatrix} n \\ k \end{bmatrix} - q^{n-k+1} \begin{bmatrix} n \\ k-1 \end{bmatrix} = \begin{bmatrix} n \\ k \end{bmatrix} - \begin{bmatrix} n \\ k-1 \end{bmatrix} + (q^k - 1)\begin{bmatrix} n \\ k \end{bmatrix} - (q^{n-k+1} - 1)\begin{bmatrix} n \\ n-k+1 \end{bmatrix}$$

$$= \begin{bmatrix} n \\ k \end{bmatrix} - \begin{bmatrix} n \\ k-1 \end{bmatrix} + (q^n - 1)\begin{bmatrix} n-1 \\ k-1 \end{bmatrix} - (q^n - 1)\begin{bmatrix} n-1 \\ n-k \end{bmatrix} = \begin{bmatrix} n \\ k \end{bmatrix} - \begin{bmatrix} n \\ k-1 \end{bmatrix}$$

we get

$$\sum_{k=0}^{\lfloor n/2 \rfloor} \begin{bmatrix} n \\ k \end{bmatrix} L_{n-2k}(x,q) = \sum \begin{bmatrix} n \\ k \end{bmatrix} F_{n-2k}(x,q) - \sum \begin{bmatrix} n \\ k \end{bmatrix} q^{n-1-2k} F_{n-2k-2}(x,q)$$

$$= \sum \left(\begin{bmatrix} n \\ k \end{bmatrix} - q^{n-2k+1}\begin{bmatrix} n \\ k-1 \end{bmatrix}\right) F_{n-2k}(x,q) = \sum_{k=0}^{\lfloor n/2 \rfloor} \frac{1}{q^k}\left(\begin{bmatrix} n \\ k \end{bmatrix} - \begin{bmatrix} n \\ k-1 \end{bmatrix}\right) F_{n-2k}(x,q) = x^n$$

and thus also (2.5).

**Remark**

The simplest recursion for $F_n(x,q)$ for a fixed number $x$ is (cf. [2])

$$F_n(x,q) = xF_{n-1}(x,q) - q^{n-2}xF_{n-3}(x,q) + q^{n-3}F_{n-4}(x,q). \tag{2.9}$$

Let us also note that comparing coefficients gives the recursion

$$F_n(x,q) = xF_{n-1}(x,q) - (\sqrt{q})^{n-2} F_{n-2}(\sqrt{q}x, q). \tag{2.10}$$

We will also consider the polynomials

$$P_n(x,q) = F_{2n}(\sqrt{x}, q) = \sum_{k=0}^{n}(-1)^{n-k} q^{\binom{n-k}{2}} \begin{bmatrix} n+k \\ 2k \end{bmatrix} x^k \tag{2.11}$$

and

$$Q_n(x,q) = \frac{F_{2n+1}(\sqrt{x}, q)}{\sqrt{x}} = \sum_{k=0}^{n}(-1)^{n-k} q^{\binom{n-k}{2}} \begin{bmatrix} n+k+1 \\ 2k+1 \end{bmatrix} x^k. \tag{2.12}$$

For these polynomials we get

$$x^n = \sum_{k=0}^{n} \frac{1}{q^{n-k}}\left(\begin{bmatrix} 2n \\ n-k \end{bmatrix} - \begin{bmatrix} 2n \\ n-k-1 \end{bmatrix}\right) P_k(x,q) \tag{2.13}$$

and

$$x^n = \sum_{k=0}^{n} \frac{1}{q^{n-k}}\left(\begin{bmatrix} 2n+1 \\ n-k \end{bmatrix} - \begin{bmatrix} 2n+1 \\ n-k-1 \end{bmatrix}\right) Q_k(x,q). \tag{2.14}$$

If we define linear functionals $L_0$ and $L_1$ by

$$L_0(P_n(x,q)) = [n=0],$$
$$L_1(Q_n(x,q)) = [n=0],$$
(2.15)

then (2.13) and (2.14) give

$$L_0(x^n) = C_n(q),$$
$$L_1(x^n) = \frac{1+q}{1+q^{n+1}} C_{n+1}(q).$$
(2.16)

Analogously let

$$R_n(x,q) = L_{2n}(\sqrt{x},q) = \sum_{k=0}^{n} (-1)^{n-k} q^{\binom{n-k}{2}} \frac{[2n]}{[n+k]} \begin{bmatrix} n+k \\ 2k \end{bmatrix} x^k$$
(2.17)

and

$$S_n(x,q) = \frac{L_{2n+1}(\sqrt{x},q)}{\sqrt{x}} = \sum_{k=0}^{n} (-1)^{n-k} q^{\binom{n-k}{2}} \frac{[2n+1]}{[n+k+1]} \begin{bmatrix} n+k+1 \\ 2k+1 \end{bmatrix} x^k.$$
(2.18)

This implies

$$x^n = \sum_{k=0}^{n} \begin{bmatrix} 2n \\ n-k \end{bmatrix} R_k(x,q),$$
$$x^n = \sum_{k=0}^{n} \begin{bmatrix} 2n+1 \\ n-k \end{bmatrix} S_k(x,q).$$
(2.19)

Let $M_0(R_n(x,q)) = [n=0]$ and $M_1(S_n(x,q)) = [n=0]$. Then we get

$$M_0(x^n) = \begin{bmatrix} 2n \\ n \end{bmatrix},$$
$$M_1(x^n) = \begin{bmatrix} 2n+1 \\ n \end{bmatrix}.$$
(2.20)

By comparing coefficients we get

$$R_n(x,q) = Q_n(x,q) - q^{2n-2} Q_{n-2}(x,q).$$
(2.21)

## 3. q-Narayana polynomials as moments

In the following we extend the above results by introducing a new parameter $t$ as in part I.

In [4] we have defined $F_n(x,t)$ by the recursions

$$F_{2n}(x,t) = xF_{2n-1}(x,t) - F_{2n-2}(x,t),$$
$$F_{2n+1}(x,t) = xF_{2n}(x,t) - tF_{2n-1}(x,t) \quad (3.1)$$

and initial values $F_0(x,t) = 1$ and $F_1(x,t) = x$.

If $L$ denotes the linear functional defined by $L(F_n(x,t)) = [n=0]$, then we have

$$L(x^{2n}) = C_n(t) \text{ and } L(x^{2n+1}) = 0 \quad (3.2)$$

where

$$C_n(t) = \sum_{k \geq 0} \binom{n}{k}\binom{n}{k+1}\frac{1}{n}t^k \quad (3.3)$$

for $n > 0$ and $C_0(t) = 1$ is a Narayana polynomial.

**Theorem 1**

Let $\tau_{2n}(t) = 1$ and $\tau_{2n+1}(t) = t$ and define $F_n(x,t,q)$ by the recursion

$$F_n(x,t,q) = xF_{n-1}(x,qt,q) - q^{\left\lfloor \frac{n-1}{2} \right\rfloor}\tau_{n-2}(t)F_{n-2}(x,t,q) \quad (3.4)$$

with initial values $F_0(x,t,q) = 1$ and $F_1(x,t,q) = x$. These polynomials are explicitly given by

$$F_{2n}(x,t,q) = \sum_{k=0}^{n}(-1)^k q^{\binom{k}{2}} \sum_{j=0}^{k} \begin{bmatrix} n-j \\ k-j \end{bmatrix}\begin{bmatrix} n-k+j-1 \\ j \end{bmatrix} q^{(n-k+1)j}t^j x^{2n-2k},$$

$$F_{2n+1}(x,t,q) = \sum_{k=0}^{n}(-1)^k q^{\binom{k}{2}} \sum_{j=0}^{k} \begin{bmatrix} n-j \\ k-j \end{bmatrix}\begin{bmatrix} n-k+j \\ j \end{bmatrix} q^{(n-k+1)j}t^j x^{2n+1-2k}. \quad (3.5)$$

If $L$ denotes the linear functional defined by $L(F_n(x,t,q)) = [n=0]$, then we have

$$L(x^{2n}) = C_n(t,q) \text{ and } L(x^{2n+1}) = 0 \quad (3.6)$$

where

$$C_n(t,q) = \frac{1}{[n]}\sum_{j=0}^{n} q^{j^2+j}\begin{bmatrix} n \\ j \end{bmatrix}\begin{bmatrix} n \\ j+1 \end{bmatrix}t^j \quad (3.7)$$

for $n > 0$ is a $q$-Narayana polynomial.

**Proof**

To prove (3.5) we have only to show that (3.4) holds. This follows by comparing coefficients.

For $t=1$ we get by using $\begin{bmatrix}-a\\k\end{bmatrix}=\begin{bmatrix}a+k-1\\k\end{bmatrix}(-q^a)^k q^{-\binom{k}{2}}$ and the $q-$Vandermonde formula

$$\sum_{j=0}^{k}\begin{bmatrix}n-j\\k-j\end{bmatrix}\begin{bmatrix}n-k+j-1\\j\end{bmatrix}q^{(n-k+1)j}$$

$$=\sum_{j=0}^{k}\begin{bmatrix}j-n+k-j-1\\k-j\end{bmatrix}(-q^{n-j})^{k-j}q^{-\binom{k-j}{2}}\begin{bmatrix}-n+k-j+1+j-1\\j\end{bmatrix}q^{-\binom{j}{2}}(-q^{n-k+j-1})^j q^{(n-k+1)j}$$

$$=(-1)^k q^{2nk-k^2-\binom{k}{2}}\sum_{j=0}^{k}\begin{bmatrix}-n+k-1\\k-j\end{bmatrix}\begin{bmatrix}-n+k\\j\end{bmatrix}q^{(k-j)(k-n-j)}=\begin{bmatrix}-2n+2k-1\\k\end{bmatrix}(-1)^k q^{2nk-k^2-\binom{k}{2}}=\begin{bmatrix}2n-k\\k\end{bmatrix}$$

and analogously

$$\sum_{j=0}^{k}\begin{bmatrix}n-j\\k-j\end{bmatrix}\begin{bmatrix}n-k+j\\j\end{bmatrix}q^{(n-k+1)j}=\begin{bmatrix}2n+1-k\\k\end{bmatrix}.$$

This implies that

$$F_n(x,1,q)=F_n(x,q). \tag{3.8}$$

To prove (3.6) let $a(n,k,t,q)$ be the uniquely defined numbers such that

$$\sum_{k=0}^{n}a(n,k,t,q)F_k(x,t,q)=x^n. \tag{3.9}$$

By (3.4) we get

$$a(n,k,t,q)=a(n-1,k-1,qt,q)+q^{\lfloor\frac{k+1}{2}\rfloor}\tau_k(t)a(n-1,k+1,qt,q) \tag{3.10}$$

with initial values $a(n,-1,t,q)=0$ and $a(0,k,t,q)=[k=0]$.

This implies that

$$a(2n+1,2k+1,t,q)=\frac{1}{q^{n-k}}\sum_{j=0}^{n-k}q^{j^2+(k+1)j}\left(\begin{bmatrix}n\\j\end{bmatrix}\begin{bmatrix}n+1\\j+k+1\end{bmatrix}-\begin{bmatrix}n+1\\j\end{bmatrix}\begin{bmatrix}n\\j+k+1\end{bmatrix}\right)t^j,$$

$$a(2n,2k,t,q)=\frac{1}{q^{n-k}}\sum_{j=0}^{n-k}q^{j^2+(k+1)j}\left(\begin{bmatrix}n-1\\j\end{bmatrix}\begin{bmatrix}n+1\\j+k+1\end{bmatrix}-\begin{bmatrix}n\\j\end{bmatrix}\begin{bmatrix}n\\j+k+1\end{bmatrix}\right)t^j \tag{3.11}$$

and $a(n,k,t,q)=0$ else.

In order to show this we must verify that

$$a(2n,2k,t,q)=a(2n-1,2k-1,qt,q)+q^k a(2n-1,2k+1,qt,q) \tag{3.12}$$

and

$$a(2n+1,2k+1,t,q)=a(2n,2k,qt,q)+q^{k+1}ta(2n,2k+2,qt,q). \tag{3.13}$$

(3.12) follows from

$$a(2n-1, 2k-1, qt, q) + q^k a(2n-1, 2k+1, qt, q) = \frac{1}{q^{n-k}} \sum_{j=0}^{n-k} q^{j^2+kj} \left( \begin{bmatrix} n-1 \\ j \end{bmatrix} \begin{bmatrix} n \\ j+k \end{bmatrix} - \begin{bmatrix} n \\ j \end{bmatrix} \begin{bmatrix} n-1 \\ j+k \end{bmatrix} \right) q^j t^j$$

$$+ \frac{q^k}{q^{n-k-1}} \sum_{j=0}^{n-k} q^{j^2+kj+j} \left( \begin{bmatrix} n-1 \\ j \end{bmatrix} \begin{bmatrix} n \\ j+k+1 \end{bmatrix} - \begin{bmatrix} n \\ j \end{bmatrix} \begin{bmatrix} n-1 \\ j+k+1 \end{bmatrix} \right) q^j t^j$$

$$= \frac{1}{q^{n-k}} \sum_{j=0}^{n-k} q^{j^2+kj+j} \left( \begin{bmatrix} n-1 \\ j \end{bmatrix} \begin{bmatrix} n+1 \\ j+k+1 \end{bmatrix} - \begin{bmatrix} n \\ j \end{bmatrix} \begin{bmatrix} n \\ j+k+1 \end{bmatrix} \right) t^j = a(2n, 2k, t, q).$$

(3.13) follows from

$$a(2n, 2k, qt, q) + q^{k+1} t a(2n, 2k+2, qt, q) = \frac{1}{q^{n-k}} \sum_{j=0}^{n-k} q^{j^2+(k+1)j} \left( \begin{bmatrix} n-1 \\ j \end{bmatrix} \begin{bmatrix} n+1 \\ j+k+1 \end{bmatrix} - \begin{bmatrix} n \\ j \end{bmatrix} \begin{bmatrix} n \\ j+k+1 \end{bmatrix} \right) q^j t^j$$

$$+ \frac{q^{k+1}}{q^{n-k-1}} \sum_{j=0}^{n-k} q^{j^2+(k+2)j} \left( \begin{bmatrix} n-1 \\ j \end{bmatrix} \begin{bmatrix} n+1 \\ j+k+2 \end{bmatrix} - \begin{bmatrix} n \\ j \end{bmatrix} \begin{bmatrix} n \\ j+k+2 \end{bmatrix} \right) q^j t^{j+1}$$

$$= \frac{1}{q^{n-k}} \sum_{j=0}^{n-k} q^{j^2+kj+2j} \left( \begin{bmatrix} n-1 \\ j \end{bmatrix} \begin{bmatrix} n+1 \\ j+k+1 \end{bmatrix} - \begin{bmatrix} n \\ j \end{bmatrix} \begin{bmatrix} n \\ j+k+1 \end{bmatrix} \right) t^j + \frac{1}{q^{n-k}} \sum_{j=0}^{n-k} q^{j^2+kj+j} \left( \begin{bmatrix} n-1 \\ j-1 \end{bmatrix} \begin{bmatrix} n+1 \\ j+k+1 \end{bmatrix} - \begin{bmatrix} n \\ j-1 \end{bmatrix} \begin{bmatrix} n \\ j+k+1 \end{bmatrix} \right)$$

$$= \frac{1}{q^{n-k}} \sum_{j=0}^{n-k} q^{j^2+kj+j} \left( \begin{bmatrix} n \\ j \end{bmatrix} \begin{bmatrix} n+1 \\ j+k+1 \end{bmatrix} - \begin{bmatrix} n+1 \\ j \end{bmatrix} \begin{bmatrix} n \\ j+k+1 \end{bmatrix} \right) t^j = a(2n+1, 2k+1, t, q)$$

As special case we get for $n > 0$

$$a(2n, 0, t, q) = \frac{1}{q^n} \sum_{j=0}^{n} q^{j^2+j} \left( \begin{bmatrix} n-1 \\ j \end{bmatrix} \begin{bmatrix} n+1 \\ j+1 \end{bmatrix} - \begin{bmatrix} n \\ j \end{bmatrix} \begin{bmatrix} n \\ j+1 \end{bmatrix} \right) t^j = C_n(t, q). \quad (3.14)$$

$C_n(t, q)$ is related to the $q$-Catalan numbers $c_n(\lambda; q)$ of J. Fürlinger and J. Hofbauer [5]. They have shown that

$$C_n(1, q) = C_n(q) = \frac{1}{[n+1]} \begin{bmatrix} 2n \\ n \end{bmatrix}. \quad (3.15)$$

This result follows again from Theorem 1 because of (3.8) and (2.3).

**Remark**

Let us also consider the polynomials $P_n(x, t, q) = F_{2n}(\sqrt{x}, t, q)$ and

$$Q_n(x, t, q) = \frac{F_{2n+1}(\sqrt{x}, t, q)}{\sqrt{x}}.$$

**Corollary 1.1**

*Let*

$$Q_n(x,t,q) = \sum_{k=0}^{n} (-1)^k q^{\binom{k}{2}} \sum_{j=0}^{k} \begin{bmatrix} n-j \\ k-j \end{bmatrix} \begin{bmatrix} n-k+j \\ j \end{bmatrix} (q^{n-k+1}t)^j x^{n-k} \qquad (3.16)$$

*and*

$$B_{n,k}(t,q) = \frac{1}{q^{n-k}} \sum_{j=0}^{n-k} q^{j^2+j+kj} \left( \begin{bmatrix} n \\ j \end{bmatrix} \begin{bmatrix} n+1 \\ k+j+1 \end{bmatrix} - \begin{bmatrix} n+1 \\ j \end{bmatrix} \begin{bmatrix} n \\ k+j+1 \end{bmatrix} \right) t^j. \qquad (3.17)$$

*Then*

$$\sum_{k=0}^{n} B_{n,k}(t,q) Q_k(x,t,q) = x^n \qquad (3.18)$$

*with*

$$B_{n,0}(t,q) = \frac{1}{[n+1]} \sum_{j=0}^{n} q^{j^2} \begin{bmatrix} n+1 \\ j \end{bmatrix} \begin{bmatrix} n+1 \\ j+1 \end{bmatrix} t^j = C_{n+1}\left(\frac{t}{q}, q\right). \qquad (3.19)$$

Note that by [5] $C_n\left(\frac{1}{q}, q\right) = \frac{1+q}{1+q^n} C_n(q)$ in accordance with (2.16).

$B_{n,k}(t,q)$ can also be written as $B_{n,k}(t,q) = \frac{[k+1]}{[n+1]} \sum_{j=0}^{n-k} q^{j^2+kj} \begin{bmatrix} n+1 \\ j \end{bmatrix} \begin{bmatrix} n+1 \\ k+j+1 \end{bmatrix} t^j.$

For $t=1$ $B_{n,k}(t,q)$ reduces to $B_{n,k}(1,q) = \frac{1}{q^{n-k}} \left( \begin{bmatrix} 2n+1 \\ n-k \end{bmatrix} - \begin{bmatrix} 2n+1 \\ n-k-1 \end{bmatrix} \right).$

**Corollary 1.2**

*Let*

$$P_n(x,t,q) = \sum_{k=0}^{n} (-1)^{n-k} q^{\binom{n-k}{2}} \sum_{j=0}^{n} \begin{bmatrix} n-j \\ k \end{bmatrix} \begin{bmatrix} k+j-1 \\ j \end{bmatrix} q^{j(k+1)} t^j x^k \qquad (3.20)$$

*and*

$$A_{n,k}(t,q) = \frac{1}{q^{n-k}} \sum_{j=0}^{n-k} q^{j^2+(k+1)j} \left( \begin{bmatrix} n-1 \\ j \end{bmatrix} \begin{bmatrix} n+1 \\ j+k+1 \end{bmatrix} - \begin{bmatrix} n \\ j \end{bmatrix} \begin{bmatrix} n \\ j+k+1 \end{bmatrix} \right) t^j. \qquad (3.21)$$

*Then*

$$\sum_{k=0}^{n} A_{n,k}(t,q) P_k(x,t,q) = x^n \qquad (3.22)$$

*with*

$$A_{n,0}(t,q) = \sum_{k=0}^{n} \begin{bmatrix} n \\ j \end{bmatrix} \begin{bmatrix} n \\ j+1 \end{bmatrix} \frac{1}{[n]} q^{j^2+j} t^j = C_n(t,q). \tag{3.23}$$

The first terms of $A_{n,k}(t,q)$ are

```
1
1                        1
1 + q² t                  1 + q + q² t                                          1
1 + q² t + q³ t + q⁴ t + q⁶ t²   1 + q + q² + q² t + q³ t + 2 q⁴ t + q⁵ t + q⁶ t²   1 + q + q² + q³ t + q⁴ t   1
```

## 2.1. q-Narayana polynomials of type B

In [4] we have seen that the orthogonal polynomials $L_n(x,t)$ whose moments are the Narayana polynomials $M_n(t) = \sum_{k=0}^{n} \binom{n}{k}^2 t^k$ of type B, satisfy the recurrence

$$L_n(x,t) = xL_{n-1}(x,t) - \tau_{n-2}(t)L_{n-2}(x,t) \tag{3.24}$$

with initial values $L_0(x,t) = 1$ and $L_1(x,t) = x.$

Here we have

$\tau_0(t) = 1+t,$

$\tau_{2n}(t) = \dfrac{1+t^{n+1}}{1+t^n}$ for $n > 0,$

$\tau_{2n+1}(t) = \dfrac{t(1+t^n)}{1+t^{n+1}}.$

We now show that there exists a natural $q$ − analogue of $L_n(x,t,q)$ with $L_n(x,t,1) = L_n(x,t)$ and which for $t = 1$ reduces to the $q$ − Lucas polynomials:

**Theorem 2**

*Let*

$$\tau_0(t,q) = 1 + qt,$$
$$\tau_{2n}(t,q) = q^n \frac{1+q^{n+1}t^{n+1}}{1+q^n t^n}, \tag{3.25}$$
$$\tau_{2n+1}(t,q) = q^{n+1} t \frac{1+t^n}{1+t^{n+1}}$$

*and define $L_n(x,t,q)$ by*

$$L_n(x,t,q) = xL_{n-1}(x,qt,q) - \tau_{n-2}(t,q)L_{n-2}(x,t,q) \tag{3.26}$$

with $L_0(x,t,q) = 1$ and $L_1(x,t,q) = x$.

Let $M$ denote the linear functional defined by $M(L_n(x,t,q)) = [n=0]$, then we get

$$M(x^{2n}) = M_n(t,q) = \sum_{j=0}^{n} q^{j^2} \begin{bmatrix} n \\ j \end{bmatrix}^2 t^j. \tag{3.27}$$

**Proof**

Let $a(n,k,t,q)$ satisfy

$$a(n,k,t,q) = a(n-1,k-1,qt,q) + \tau_k(t,q)a(n-1,k+1,qt,q) \tag{3.28}$$

with $a(n,-1,t,q) = 0$ and $a(0,k,t,q) = [k=0]$.

Then (3.26) implies that

$$\sum_{k=0}^{n} a(n,k,t,q)L_k(x,t,q) = x^n. \tag{3.29}$$

By induction it is easy to verify that

$$a(2n,2k,t,q) = \sum_{j=0}^{n-k} q^{j(j+k)} \begin{bmatrix} n \\ j \end{bmatrix} \begin{bmatrix} n \\ j+k \end{bmatrix} t^j,$$

$$a(2n+1,2k+1,t,q) = \frac{1}{1+t^{k+1}} \sum_{j=0}^{n-k} \begin{bmatrix} n \\ k+j \end{bmatrix} \begin{bmatrix} n+1 \\ j \end{bmatrix} \left( q^{j(j+k)} t^j + q^{(n+1-j)(n-j-k)} t^{n+1-j} \right) \tag{3.30}$$

$$= \frac{1}{1+t^{k+1}} \left( \sum_{j=0}^{n-k} \begin{bmatrix} n \\ k+j \end{bmatrix} \begin{bmatrix} n+1 \\ j \end{bmatrix} q^{j(j+k)} t^j + \sum_{j=k+1}^{n+1} \begin{bmatrix} n \\ j-k-1 \end{bmatrix} \begin{bmatrix} n+1 \\ j \end{bmatrix} q^{j(j-k-1)} t^j \right)$$

and $a(n,k,t,q) = 0$ else.

If $M(L_n(x,t,q)) = [n=0]$ then (3.29) implies

$$M(x^{2n}) = a(2n,0,t,q) = \sum_{j=0}^{n} q^{j^2} \begin{bmatrix} n \\ j \end{bmatrix}^2 t^j = M_n(t,q).$$

By $q$-Vandermonde we see that $a(2n,2k,1,q) = \begin{bmatrix} 2n \\ n-k \end{bmatrix}$ and $a(2n+1,2k+1,1,q) = \begin{bmatrix} 2n+1 \\ n-k \end{bmatrix}$.

Comparing with (2.19) we conclude that $L_n(x,1,q)$ is a $q$-Lucas polynomial. Further it is clear that $L_n(x,t,1) = L_n(x,t)$.

Let $R_n(x,t,q) = L_{2n}(\sqrt{x},t,q)$.

In order to get a formula for $R_n(x,t,q)$ observe that

$$B_{n,k}(t,q) = D_{n,k}(t,q) - q^{2k+2}tD_{n,k+2}(t,q), \tag{3.31}$$

which is equivalent with the easily verified identity

$$\frac{[k+1]}{[n+1]}\begin{bmatrix}n+1\\j\end{bmatrix}\begin{bmatrix}n+1\\k+j+1\end{bmatrix} = \begin{bmatrix}n\\j\end{bmatrix}\begin{bmatrix}n\\k+j\end{bmatrix} - q^{k+1}\begin{bmatrix}n\\j-1\end{bmatrix}\begin{bmatrix}n\\k+j+1\end{bmatrix}.$$

By (3.18) this implies

$$\sum_{k=0}^{n} D_{n,k}(t,q)\left(Q_k(x,t,q) - q^{2k-2}tQ_{k-2}(x,t,q)\right)$$

$$= \sum_{k=0}^{n} D_{n,k}(t,q)Q_k(x,t,q) - \sum_{k=0}^{n} q^{2k+2}tD_{n,k+2}(t,q)Q_k(x,t,q) = \sum_{k=0}^{n} B_{n,k}(t,q)Q_k(x,t,q) = x^n.$$

Comparing with (3.37) we see that

$$R_n(x,t,q) = Q_n(x,t,q) - q^{2n-2}tQ_{n-2}(x,t,q) \tag{3.32}$$

if we let $Q_{-2}(x,t,q) = Q_{-1}(x,t,q) = 0$.

The first terms of $(R_n(x,t,q))_{n\geq 0}$ are

```
{1, -1 - q t + x, q + q³ t² - x - q x - q² t x - q³ t x + x²,
-q³ - q⁶ t³ + q x + q² x + q³ x + q³ t x + q⁴ t x + q⁵ t x + q⁵ t² x +
q⁶ t² x + q⁷ t² x - x² - q x² - q² x² - q³ t x² - q⁴ t x² - q⁵ t x² + x³}
```

From (3.32) and (3.16) we get the formula

$$R_n(x,t,q) = \sum_{k=0}^{n} (-1)^k q^{\binom{k}{2}} \begin{bmatrix}n\\k\end{bmatrix} c(n,k,t) x^{n-k} \tag{3.33}$$

with

$$c(n,k,t) = \sum_{j=0}^{k} \begin{bmatrix}k\\j\end{bmatrix} q^{(n+1-k)j} \frac{\begin{bmatrix}n+j-k-1\\j\end{bmatrix}}{\begin{bmatrix}n-1\\j\end{bmatrix}} t^j \quad \text{for } k < n, \tag{3.34}$$

$$c(n,n,t) = 1 + q^n t^n.$$

This can also be written as

$$R_k(x,t,q) = (-1)^k q^{\binom{k}{2}}\left(1+q^k t^k\right) + \sum_{\ell=1}^{k} (-1)^{k-\ell} q^{\binom{k-\ell}{2}} \begin{bmatrix}k\\\ell\end{bmatrix} x^\ell \sum_{j=0}^{k-\ell} q^{(\ell+1)j} \begin{bmatrix}k-\ell\\j\end{bmatrix} \frac{\begin{bmatrix}\ell+j-1\\j\end{bmatrix}}{\begin{bmatrix}k-1\\j\end{bmatrix}} t^j. \tag{3.35}$$

By (3.29) and (3.30) we get

**Corollary 2.1**

Let

$$D_{n,k}(t,q) = a(2n,2k,t,q) = \sum_{j=0}^{n-k} q^{j(j+k)} \begin{bmatrix} n \\ j \end{bmatrix}\begin{bmatrix} n \\ k+j \end{bmatrix} t^j. \qquad (3.36)$$

Then

$$\sum_{k=0}^{n} D_{n,k}(t,q) R_k(x,t,q) = x^n. \qquad (3.37)$$

Let $M_0$ be the linear functional defined by $M_0(R_n(x,t,q)) = [n=0]$. Then

$$M_0(x^n) = D_{n,0}(t,q) = B_n(t,q) = \sum_{j=0}^{n} q^{j^2} \begin{bmatrix} n \\ j \end{bmatrix}^2 t^j. \qquad (3.38)$$

Let now $S_n(x,t,q) = \dfrac{L_{2n+1}(\sqrt{x},t,q)}{\sqrt{x}}.$

From (3.26) we get for $n>0$

$$S_n(x,t,q) = \frac{1}{x(1+t^n)}\left((1+t^n)R_{n+1}\left(x,\frac{t}{q},q\right) + q^n(1+t^{n+1})R_n\left(x,\frac{t}{q},q\right)\right). \qquad (3.39)$$

**Corollary 2.2**

Let

$$E_{n,k}(t,q) = a(2n+1,2k+1,t,q) = \frac{1}{1+t^{k+1}} \sum_{j=0}^{n-k} \begin{bmatrix} n \\ k+j \end{bmatrix}\begin{bmatrix} n+1 \\ j \end{bmatrix}\left(q^{j(j+k)}t^j + q^{(n+1-j)(n-j-k)}t^{n+1-j}\right).$$
$$(3.40)$$

Then

$$\sum_{k=0}^{n} E_{n,k}(t,q) S_k(x,t,q) = x^n. \qquad (3.41)$$

Let $M_1$ be the linear functional defined by $M_1(S_n(x,t,q)) = [n=0]$. Then

$$M_1(x^n) = E_{n,0}(t,q) = \frac{1}{1+t}\sum_{j=0}^{n+1} q^{j^2-j}\begin{bmatrix} n+1 \\ j \end{bmatrix}^2 t^j = \frac{M_{n+1}\left(\frac{t}{q},q\right)}{1+t}. \qquad (3.42)$$

**Remark**

For the numbers $D_{n,k}(t,q)$ there exists an analogue of the Catalan-Stieltjes matrix for orthogonal polynomials:

$$D_{n,0}(t,q) = (1+q^n t)D_{n-1,0}(qt,q) + q(1+q^n)tD_{n-1,1}(qt,q) \tag{3.43}$$

and

$$D_{n,k}(t,q) = D_{n-1,k-1}(qt,q) + q^k(1+q^n t)D_{n-1,k}(qt,q) + q^{n+2k+1}tD_{n-1,k+1}(qt,q). \tag{3.44}$$

Let us mention some curious conjectures: Let $\Delta_{q,t} f(t) = \dfrac{f(t)-f(qt)}{(1-q)t}$ be the $q$-differential operator with respect to the variable $t$. Then

$$\sum_{k=0}^{n} \Delta_{q,t}^m \left(D_{n,k}(x,t)\right) R_k(x, q^m t, q) = q^{m^2}[m]! \begin{bmatrix} n \\ m \end{bmatrix} \sum_{j=0}^{n-m} q^{mj} c_{n-m-j}\left(q^{2j}t, m, q\right) x^j. \tag{3.45}$$

Here $c_n(t,m,q) = 0$ for $n < 0$, $c_0(t,m,q) = 1$ and for $n > 0$

$$c_n(t,m,q) = \sum_{k=0}^{n-1} \begin{bmatrix} n-1 \\ k \end{bmatrix} \begin{bmatrix} n+m \\ k+m \end{bmatrix} \frac{[m]}{[n+m]} q^{2km+k^2+k} t^k. \tag{3.46}$$

For $m = 1$ we get

$$c_n(t,1,q) = \sum_{k=0}^{n-2} \begin{bmatrix} n-2 \\ k \end{bmatrix} \begin{bmatrix} n \\ k+1 \end{bmatrix} \frac{1}{[n]} q^{2k+k^2+k} t^k = C_{n-1}(q^2 t)$$

Thus

$$\sum_{k=0}^{n} \Delta_{q,t}\left(D_{n,k}(x,t)\right) R_k(x, qt, q) = \begin{bmatrix} n \\ 1 \end{bmatrix} \sum_{j=0}^{n-1} q^{j+1} C_{n-j-1}\left(q^{2j+2}t, q\right) x^j. \tag{3.47}$$

It should be noted that the numbers $c_n(t,m,1) = \sum_{k=0}^{n-m-1} \binom{n-1}{k}\binom{n+m}{k+m}\dfrac{m}{n+m} t^k$

are the coefficients of the powers $C(x,t)^m$ of the generating function

$C(x,t) = \sum_{n \geq 0} C_n(t) x^n$ of the Narayana polynomials (cf. [4]).

Other such identities are

$$\sum_{k=0}^{n} \Delta_{q,t} A_{n,k}(t,q) P_k(x, qt, q) = \sum_{j=2}^{n} q^j \begin{bmatrix} j-1 \\ 1 \end{bmatrix} C_{n-j}\left(q^{2j}t, q\right) x^{j-1}. \tag{3.48}$$

$$\sum_{k=0}^{n} \Delta_{q,t} B_{n,k}(t,q) Q_k(x, qt, q) = \sum_{j=1}^{n} q^j \begin{bmatrix} j \\ 1 \end{bmatrix} C_{n-j}\left(q^{2j+1}t, q\right) x^{j-1}. \tag{3.49}$$

$$\sum_{k=0}^{n} \Delta_{q,t}^{m}\left(A_{n,k}(x,t)\right)P_k(x,q^m t,q) = q^{m^2+m}\prod_{j=1}^{m-1}[n-j]\sum_{j=0}^{n-m}q^{mj}[j+1]c_{n-m-1-j}\left(q^{2j+2}t,m,q\right)x^{j+1}.$$
(3.50)

$$\sum_{k=0}^{n} \Delta_{q,t}^{m}\left(B_{n,k}(x,t)\right)Q_k(x,q^m t,q) = q^{m^2}\prod_{j=1}^{m-1}[n+1-j]\sum_{j=0}^{n-m}q^{mj}[j+1]c_{n-m-j}\left(q^{2j+1}t,m,q\right)x^{j}.$$
(3.51)